\documentclass[12pt,leqno]{article}

\usepackage{latexsym,amsmath,amssymb,amsthm,amsfonts}
\usepackage[active]{srcltx}     
\usepackage[all]{xy}

\renewcommand{\epsilon}{\varepsilon}

\newcommand{\R}{\mathbb{R}}

\newcommand{\p}{\mathcal{P}}

\newtheorem{theorem}{Theorem}
\newtheorem{lemma}{Lemma}

\newtheorem{corollary}{Corollary}

\begin{document}

\author{Andriy Bondarenko,
Danylo Radchenko, and Maryna Viazovska}
\title{Optimal asymptotic bounds for spherical designs}
\date{}
\maketitle
\begin{abstract}
In this paper we prove the conjecture of Korevaar and Meyers: for
each $N\ge c_dt^d$ there exists a spherical $t$-design in the sphere
$S^d$ consisting of $N$ points, where $c_d$ is a constant depending
only on~$d$.
\end{abstract}

\newpage
\section{Introduction}

Let $S^d$ be the unit sphere in $\R^{{d+1}}$ with the Lebesgue
measure $\mu_d$ normalized by $\mu_d(S^d)=1$.

A set of points $x_1,\ldots,x_N\in S^d$ is called a {\it spherical
$t$-design} if
$$ \int_{S^d}P(x)\,d\mu_d(x)=\frac 1N\sum_{i=1}^N P(x_i) $$ for all
algebraic polynomials in $d+1$ variables, of total degree at most
$t$. The concept of a spherical design was introduced by Delsarte,
Goethals, and Seidel~\cite{DGS}. For each $t, d\in {\mathbb N}$
denote by $N(d,t)$ the minimal number of points in a spherical
$t$-design in $S^d$. The following lower bound
\begin{equation}
\label{hh} N(d,t)\ge \begin{cases}\displaystyle {{d+k}\choose{d}}+{{d+k-1}\choose{d}}&\text{if $
t=2k$,}\\&\\ \displaystyle   2\,{{d+k}\choose{d}}& \text{if $ t=2k+1$,}\end{cases}
\end{equation}
is proved in \cite{DGS}.

Spherical $t$-designs attaining this bound are called tight. The
vertices of a regular $t+1$-gon form a tight spherical $t$-design in
the circle, so $N(1,t)=t+1$. Exactly eight tight spherical designs
are known for $d\ge 2$ and $t\ge 4$. All such configurations of
points are highly symmetrical, and optimal from many different
points of view (see Cohn, Kumar~\cite{CK} and Conway,
Sloane~\cite{CS}). Unfortunately, tight designs rarely exist. In
particular, Bannai and Damerell~\cite{Bannai1, Bannai2} have shown
that tight spherical designs with $d\ge 2$ and $t\ge 4$ may exist
only for $t=4$, $5$, $7$ or $11$. Moreover, the only tight
$11$-design is formed by minimal vectors of the Leech lattice in
dimension $24$. The bound~\eqref{hh} has been improved by Delsarte's
linear programming method for most pairs $(d,\,t)$; see~\cite{Yu}.

On the other hand, Seymour and Zaslavsky~\cite{SZ} have proved that
spherical $t$-designs exist for all $d$, $t\in {\mathbb N}$.
However, this proof is nonconstructive and gives no idea of how big
$N(d,t)$ is. So, a natural question is to ask how $N(d,t)$ differs
from the tight bound~\eqref{hh}. Generally, to find the exact value
of $N(d,t)$ even for small $d$ and $t$ is a surprisingly hard
problem. For example, everybody believes that 24 minimal vectors of
the $D_4$ root lattice form a $5$-design with minimal number of
points in $S^3$, although it is only proved that $22\le N(3,5)\le
24$; see~\cite{BDN}. Further, Cohn, Conway, Elkies, and
Kumar~\cite{CCEK} conjectured that every spherical $5$-design
consisting of $24$ points in $S^3$ is in a certain $3$-parametric
family. Recently, Musin~\cite{M} has solved a long standing problem
related to this conjecture. Namely, he proved that the kissing
number in dimension $4$ is $24$.

In this paper we focus on asymptotic upper bounds on $N(d,t)$ for
fixed $d\ge 2$ and $t\to\infty$. Let us give a brief history of this
question. First, Wagner~\cite{Wag} and Bajnok~\cite{B} proved that
$N(d,t)\le C_dt^{Cd^4}$ and $N(d,t)\le C_dt^{Cd^3}$, respectively.
Then, Korevaar and Meyers \cite{KM} have improved these inequalities
by showing that $N(d,t)\le C_dt^{(d^2+d)/2}$. They have also
conjectured that
$$
 N(d,t)\le C_dt^d.
$$
Note that \eqref{hh} implies $N(d,t)\ge c_dt^d$. Here and in what
follows we denote by $C_d$ and $c_d$ sufficiently large and
sufficiently small positive constants depending only on $d$,
respectively.

The conjecture of Korevaar and Meyers attracted the interest of many
mathematicians. For instance, Kuijlaars and Saff~\cite{SK}
emphasized the importance of this conjecture for $d=2$, and revealed
its relation to minimal energy problems. Mhaskar, Narcowich, and
Ward~\cite{MNW} have constructed positive quadrature formulas in
$S^d$ with $C_dt^d$ points having {\em almost} equal weights. Very
recently, Chen, Frommer, Lang, Sloan, and Womersley~\cite{CFL,CW}
gave a computer-assisted proof that spherical $t$-designs with
$(t+1)^2$ points exist in $S^2$ for $t\le 100$.

For $d=2$, there is an even stronger conjecture by Hardin and
Sloane~\cite{HS} saying that $N(2,t)\le\frac12t^2+o(t^2)$ as
$t\to\infty$. Numerical evidence supporting the conjecture was also
given.

In~\cite{BV1}, we have suggested a nonconstructive approach for
obtaining asymptotic bounds for $N(d,t)$ based on the application of
the Brouwer fixed point theorem. This led to the following result:
\\{\it For each $N\ge C_dt^\frac{2d(d+1)}{d+2}$ there
exists a spherical $t$-design in $S^d$ consisting of $N$ points.}\\
Instead of the Brouwer fixed point theorem we use in this paper the
following result from the Brouwer degree theory~\cite[Th. 1.2.6, Th. 1.2.9]{OCC}.\\
{\sc Theorem A. }{\it Let $f: \R^n\to \R^n$ be a continuous mapping
and $\Omega$ an open bounded subset, with boundary $\partial\Omega$,
such that  $0\in\Omega\subset \R^n$. If $(x, f(x))> 0$ for all
$x\in\partial \Omega$, then there exists $x\in \Omega$ satisfying $f(x)=0$.}\\
We employ this theorem to prove the conjecture of Korevaar and
Meyers.
\begin{theorem}\label{main}
For each $N\ge C_dt^d$ there exists a spherical $t$-design in $S^d$
consisting of $N$ points.
\end{theorem}
Note that Theorem 1 is slightly stronger than the original
conjecture because it guarantees the existence of spherical
$t$-designs for {\it each} $N$ greater than $C_dt^d$.

This paper is organized as follows. In Section 2 we explain the main
idea of the proof. Then in Section 3 we present some auxiliary
results. Finally, we prove Theorem 1 in Section 4.

\section{Preliminaries and the main idea}
Let $\mathcal{P}_t$ be the Hilbert space of polynomials $P$ on
$S^{d}$ of degree at most $t$ such that $$
\int_{S^d}P(x)d\mu_d(x)=0,
$$
equipped with the usual inner product
$$
(P,Q)=\int_{S^d}P(x)Q(x)d\mu_d(x).
$$
By the Riesz representation theorem, for each point $x\in S^d$ there
exists a unique polynomial $G_x\in \p_t$ such that
$$(G_x,Q)=Q(x) \;\;\mbox{for all}\;\;Q\in\p_t.$$
Then a set of points $x_1,\ldots,x_N\in S^d$ forms a spherical
$t$-design if and only if
\begin{equation}
\label{sph} G_{x_1}+\cdots+G_{x_N}=0.
\end{equation}
For a  differentiable function $f: \R^{d+1}\to \R$ denote by
$$ \frac{\partial f}{\partial x}(x_0):=\left(\frac{\partial
f}{\partial \xi_1}(x_0),\ldots, \frac{\partial f}{\partial
\xi_{d+1}}(x_0)\right) $$ the gradient of $f$  at the point $x_0\in
\R^{d+1}$.

For a polynomial $Q\in\p_t$ we define the spherical gradient as
follows:
\begin{equation}
\label{grad} \nabla Q(x):=\frac{\partial}{\partial x}Q\left(\frac
{x}{|x|}\right),
\end{equation}
where $|\cdot|$ is the Euclidean norm in $\R^{d+1}$.

 We apply Theorem A to the open subset $\Omega$ of a vector
 space~$\p_t$,
 \begin{equation}\label{omega}\Omega:=\left\{P\in\mathcal{P}_t\left|\,\int_{S^d}|\nabla
P(x)|d\mu_d(x)<1\right.\right\}.
\end{equation}

Now we observe that the existence of a continuous mapping $F:
\mathcal{P}_t\to (S^d)^N$, such that for all $P\in\partial\Omega$
\begin{equation}\label{positive}\sum_{i=1}^N
P(x_i(P))>0,\;\mathrm{where}\;F(P)=(x_1(P),...,x_N(P)),\end{equation}
 readily implies the existence of a spherical
$t$-design in $S^d$ consisting of $N$ points. Consider a mapping
$L:(S^d)^N\to \p_t$ defined by
$$(x_1,\ldots,x_N)\ \,{\mathop{{\longrightarrow}}\limits^{L}}\ \,
G_{x_1}+\cdots+G_{x_N},$$ and the following composition mapping
$f=L\circ F: \mathcal{P}_t\to\mathcal{P}_t$. Clearly
$$
(P,f(P))=\sum_{i=1}^N P(x_i(P))
$$
for each $P\in\mathcal{P}_t$. Thus, applying Theorem A to the
mapping $f$, the vector space $ \mathcal{P}_t$, and the subset
$\Omega$ defined by \eqref{omega}, we obtain that $f(Q)=0$ for some
$Q\in\mathcal{P}_t$. Hence, by~\eqref{sph}, the components of
$F(Q)=(x_1(Q),...,x_N(Q))$
 form a spherical $t$-design in $S^d$ consisting of $N$
points.

The most naive approach to construct such $F$ is to start with a
certain well-distributed collection of points $x_i$
($i=1,\ldots,N$), put $F(0):=(x_1,\ldots,x_N)$, and then move each
point along the spherical gradient vector field of $P$. Note that
this is the most greedy way to increase each $P(x_i(P))$ and make
$\sum_{i=1}^N P(x_i(P))$ positive for each $P\in\partial\Omega$.
Following this approach we will give an explicit construction of $F$
in Section 4, which will immediately imply the proof of Theorem 1.
\section{Auxiliary results}
To construct the corresponding mapping $F$ for each $N\ge C_dt^d$ we
extensively use the following notion of an area-regular partition.

Let $\mathcal{R}=\left\{R_1,\ldots, R_N\right\}$ be a finite
collection of closed sets $R_i\subset S^d$ such that $\cup_{i=1}^N
R_i=S^d$ and $\mu_d(R_i\cap R_j)=0$ for all $1\le i<j\le N$. The
partition $\mathcal{R}$ is called area-regular if $\mu_d(R_i)=1/N$,
$i=1,\ldots,N$. The partition norm for $\mathcal{R}$ is defined by
$$ \|\mathcal{R}\|:=\max_{R\in\mathcal{R}}\mathrm{diam}\, R,$$
where $\mathrm{diam}\,R$ stands for the maximum  geodesic  distance
between two points in $R$. We need the following fact on
area-regular partitions (see Bourgain,
 Lindenstrauss~\cite{BL} and Kuijlaars, Saff~\cite{SK2}):\\
{\sc Theorem B. }{\it For each $N\in\mathbb{N}$ there exists an
area-regular partition $\mathcal{R}=\left\{R_1,\ldots, R_N\right\}$
with $\|\mathcal{R}\|\leq B_dN^{-1/d}$ for some constant $B_d$ large
enough.}

We will also use the following spherical Marcinkiewicz--Zygmund type
inequality: \\
{\sc Theorem C. }{\it There exists a constant $r_d$ such that for
each area-regular partition $\mathcal{R}=\{R_1,\ldots,R_N\}$ with
$\|\mathcal{R}\|<\frac{r_d}m$, each collection of points $x_i\in
R_i$ ($i=1,\ldots,N$), and each algebraic polynomial $P$ of total
degree $m$, the inequality
\begin{equation}
\label{Mhaskar}
\frac12\int_{S^d}|P(x)|d\mu_d(x)\le\frac1N\sum_{i=1}^N|P(x_i)|\le\frac32\int_{S^d}|P(x)|d\mu_d(x)
\end{equation}
holds.}\\
Theorem C follows naturally from the proof of Theorem 3.1
in~\cite{MNW}.
\begin{corollary} For each area-regular partition $\mathcal{R}=\{R_1,\ldots,R_N\}$
with $\|\mathcal{R}\|<\frac{r_d}{m+1}$, each collection of points
$x_i\in R_i$ ($i=1,\ldots,N$), and each algebraic polynomial $P$ of
total degree $m$,
\begin{equation}
\label{Mhaskar-mod} \frac{1}{3\sqrt{d}}\int_{S^d}|\nabla
P(x)|d\mu_d(x)\le\frac1N\sum_{i=1}^N|\nabla P(x_i)|\le
3\sqrt{d}\int_{S^d}|\nabla P(x)|d\mu_d(x).
\end{equation}
\end{corollary}
\begin{proof}
Since $|\nabla P|=\sqrt{P_1^2+\ldots+P_{d+1}^2}$ in $S^d$, where
$P_j$ are polynomials of total degree $m+1$, Corollary 1 is an
immediate consequence of~\eqref{Mhaskar} applied to $P_j$,
$j=1,\ldots,d+1$.
\end{proof}
\section{Proof of Theorem 1}
In this section we construct the map $F$ introduced in Section 2 and
thereby finish the proof of Theorem 1.

For $d, t\in{\mathbb N}$, take $C_d>(54 d B_d/r_d)^d$, where $B_d$
is as in Theorem B and $r_d$ is as in Theorem C, and fix $N\ge
C_dt^d$. Now we are in a position to give an exact construction of
the mapping $F: \mathcal{P}_{t}\to(S^d)^N$ which satisfies condition
\eqref{positive}. Take an area-regular partition
$\mathcal{R}=\left\{R_1,\ldots, R_N\right\}$ with

\begin{equation}
\label{eee}
 \|\mathcal{R}\|\le B_d N^{-1/d}<\frac{r_d}{54dt}
\end{equation}
as provided by Theorem B, and choose an arbitrary $x_i\in R_i$ for
each $i=1,\ldots, N$. Put $\epsilon=\frac{1}{6\sqrt d}$ and consider
the function
$$ h_\epsilon(u):=\begin{cases} u & \text{ if $u>\epsilon$},\\ \epsilon& \mbox{otherwise.}\end{cases} $$
Take a mapping $U:\p_t\times S^d\to{\mathbb R}^{d+1}$ such that
$$
U(P, y)=\frac{\nabla P(y)}{h_\epsilon (|\nabla P(y)|)}.
$$
For each $i=1,\ldots,N$ let $y_i:\p_t\times[0,\infty)\to S^d$ be the
map satisfying the differential equation
\begin{equation}
\label{diffur} \frac{d }{ds}y_i(P,s)=U(P,y_i(P,s))
\end{equation}
with the initial condition
$$
y_i(P,0)=x_i,
$$
for each $P\in\p_t$. Note that each mapping $y_i$ has its values in
$S^d$ by definition of spherical gradient~\eqref{grad}. Since the
mapping $U(P,y)$ is Lipschitz continuous in both $P$ and $y$, each
$y_i$ is well defined and continuous in both $P$ and $s$, where the
metric on $\p_t$ is given by the inner product.
 Finally put
\begin{equation}
\label{map} F(P)=(x_1(P),\ldots,x_N(P)):=\big(y_1(P,
\frac{r_d}{3t}),\ldots,y_N(P,\frac{r_d}{3t})\big ).
\end{equation}By definition the mapping $F$ is continuous on $\p_t$.
So, as explained in Section~2, to finish the proof of Theorem 1 it
suffices to prove
\begin{lemma}
Let $F: \mathcal{P}_{t}\to(S^d)^N$ be the mapping defined
by~\eqref{map}. Then for each $P\in\partial\Omega$,
$$
\frac1N\sum_{i=1}^N P(x_i(P))>0,
$$
where $\Omega$ is given by~\eqref{omega}.
\end{lemma}
\begin{proof}Fix $P\in\partial\Omega$. For the sake of
simplicity we write $y_i(s)$ in place of $y_i(P,s)$. By the
Newton-Leibniz formula we have
$$
\frac1N\sum_{i=1}^N P(x_i(P))=\frac
1N\sum_{i=1}^NP(y_i(r_d/3t))
$$
\begin{equation}
\label{ee}=\frac 1N\sum_{i=1}^NP(x_i) +\int_0^{r_d/3t} \frac{d}{d s}\left[\frac 1N\sum_{i=1}^N
P(y_i(s))\right]ds.
\end{equation}
Now to prove Lemma 1, we first estimate the value
$$
\left|\frac 1N\sum_{i=1}^NP(x_i)\right|
$$
from above, and then estimate the value
$$
\frac{d}{d s}\left[\frac 1N\sum_{i=1}^N P(y_i(s))\right]
$$
from below, for each $s\in [0, r_d/3t]$. We have
$$
\left|\frac
1N\sum_{i=1}^NP(x_i)\right|=\left|\sum_{i=1}^N\int_{R_i}P(x_i)-P(x)\,d\mu_d(x)\right|\le
\sum_{i=1}^N\int_{R_i}|P(x_i)-P(x)|d\mu_d(x)
$$
$$
\le\frac {\|\mathcal{R}\|}{N}\sum_{i=1}^N\max_{z\in
S^d:\,\mathrm{dist}(z,x_i)\le\|\mathcal{R}\|}|\nabla P(z)|
$$
where $\mathrm{dist}(z,x_i)$ denotes the geodesic distance between
$z$ and $x_i$. Hence, for $z_i\in S^d$ such that
$\mathrm{dist}(z_i,x_i)\le\|\mathcal{R}\|$ and
$$|\nabla P(z_i)|=\max_{z\in
S^d:\,\mathrm{dist}(z,x_i)\le\|\mathcal{R}\|}|\nabla P(z)|,$$ we
obtain
$$
\left|\frac 1N\sum_{i=1}^NP(x_i)\right|\le\frac{\|\mathcal{R}\|}{N}
\sum_{i=1}^N|\nabla P(z_i)|.
$$
Consider another area-regular partition
$\mathcal{R'}=\left\{R'_1,\ldots, R'_N\right\}$ defined by
$R'_i=R_i\cup\{z_i\}$. Clearly $\|\mathcal{R'}\|\le
2\|\mathcal{R}\|$ and so, by~\eqref{eee}, we get
$\|\mathcal{R'}\|<r_d/(27\, d\,t)$. Applying
inequality~\eqref{Mhaskar-mod} to the partition $\mathcal{R'}$ and
the collection of points $z_i$ we obtain that
\begin{equation}
\label{e1} \left|\frac 1N\sum_{i=1}^NP(x_i)\right|\le
3\sqrt{d}\,\|\mathcal{R}\|\,\int_{S^d}|\nabla P(x)|d\mu_d(x)<\frac{
r_d}{18\,\sqrt{d} \, t}
\end{equation}
for any $P\in\partial\Omega$. On the other hand, the differential
equation \eqref{diffur} implies \begin{align} \frac{d}{d
s}\left[\frac 1N\sum_{i=1}^N P(y_i(s))\right]=
 & \frac 1N\sum_{i=1}^N\frac{|\nabla P(y_i(s))|^2}{h_\epsilon (|\nabla P(y_i(s))|)
} \notag \\ \geq & \frac 1N\sum_{i:\,|\nabla P(y_i(s))|\geq \epsilon} |\nabla P(y_i(s))|\notag \\
\label{epsilon}\geq & \frac 1N\sum_{i=1}^N |\nabla
P(y_i(s))|-\epsilon.
\end{align}
Since $$\left|\frac{\nabla P(y)}{h_\epsilon (|\nabla
P(y)|)}\right|\le1$$ for each $y\in S^d$, it follows again from
\eqref{diffur} that $\left|\frac{d y_i(s)}{d s}\right|\le1$. Hence
we arrive at
$$
\mathrm{dist}(x_i,y_i(s))\leq s.
$$
Now for each $s\in [0, r_d/3t]$ consider the area-regular partition
$\mathcal{R''}=\left\{R''_1,\ldots, R''_N\right\}$ given by
$R''_i=R_i\cup\{y_i(s)\}$. By~\eqref{eee} we have
$$
\|\mathcal{R''}\|<\frac{r_d}{54dt}+\frac{r_d}{3t};
$$
so we can apply~\eqref{Mhaskar-mod} to the partition $\mathcal{R''}$
and the collection of points $y_i(s)$. This and
inequality~\eqref{epsilon} yield
$$
\frac{d}{d s}\left[\frac 1N\sum_{i=1}^N P(y_i(s))\right]\ge \frac
1N\sum_{i=1}^N |\nabla P(y_i(s))|-\frac 1{6\sqrt{d}}
$$
\begin{equation}
\label{ee2} \ge\frac 1{3\sqrt{d}}\int_{S^d}|\nabla
P(x)|d\mu_d(x)-\frac 1{6\sqrt{d}}=\frac 1{6\sqrt{d}},
\end{equation}
for each $P\in\partial\Omega$ and $s\in [0,r_d/3t]$. Finally,
equation~\eqref{ee} and inequalities \eqref{e1} and \eqref{ee2}
imply
\begin{equation}
\label{cond3} \frac1N\sum_{i=1}^N P(x_i(P))> \frac
1{6\sqrt{d}}\,\frac{r_d}{3t}-\frac{ r_d}{18\,\sqrt{d} \, t}=0.
\end{equation}
Lemma 1 is proved.
\end{proof}


{\footnotesize \noindent Centre de Recerca Matem\`atica, Campus de
Bellaterra, 08193
Bellaterra (Barcelona), Spain\\
and\\
Department of Mathematical Analysis, National Taras Shevchenko
University, str.\ Volodymyrska, 64, Kyiv, 01033, Ukraine\\
{\it Email address: andriybond@gmail.com}\\
\vspace{0.5cm}

\noindent Department of Mathematical Analysis, National Taras
Shevchenko University, str.\ Volodymyrska, 64, Kyiv, 01033, Ukraine\\
{\it Email address: danradchenko@gmail.com}\\
\vspace{0.5cm}

\noindent
Max Planck Institute for Mathematics, Vivatsgasse 7, 53111 Bonn, Germany\\
{\it Email address: viazovsk@mpim-bonn.mpg.de}}

\begin{thebibliography}{10}




\bibitem{B}
{\sc B. Bajnok}, {\em Construction of spherical $t$-designs}, Geom.
Dedicata, 43 (1992), 167--179.

\bibitem{Bannai1} {\sc E. Bannai and R.M. Damerell},{\em  Tight spherical designs I}, J. Math. Soc. Japan 31 (1979), pp.~199--207.
 \bibitem{Bannai2} {\sc E. Bannai and R.M. Damerell}, {\em Tight spherical designs. II}, J. London Math. Soc. 21 (1980), pp.~13--30.
\bibitem{BV1}
{\sc A. Bondarenko and M. Viazovska}, {\em Spherical designs via
Brouwer fixed point theorem}, SIAM J. Discrete Math., 24 (2010),
207--217.

\bibitem{BL}
{\sc J. Bourgain and J. Lindenstrauss}, {\em Distribution of points
on spheres and approximation by zonotopes}, Israel J. Math., 64
(1988), 25--31.
\bibitem{BDN}
{\sc P. Boyvalenkov, D. Danev, and S. Nikova}, {\em Nonexistence of
certain spherical designs of odd strengths and cardinalities},
Discrete Comput. Geom. 21 (1999), 143--156.

\bibitem{CCEK}
{\sc H. Cohn, J. H. Conway, N. D. Elkies, and A. Kumar}, {\em The
$D_4$ root system is not universally optimal}, Experiment. Math., 16
(2007), 313-320.

\bibitem{CK}
{\sc H. Cohn and A. Kumar}, {\em Universally optimal distribution of
points on spheres}, J. Amer. Math. Soc., 20 (2007), 99--148.

\bibitem{CFL}
{\sc X. Chen, A. Frommer, and B. Lang}, {\em Computational existence
proofs for spherical $t$-designs}, Numerische Mathematik, 117
(2011), 289--305.

\bibitem{CW}
{\sc X. Chen and R. S. Womersley}, {\em Existence of solutions to
systems of underdetermined equations and spherical designs}, SIAM
Journal on Numerical Analysis 44, 6 (2006), 2326 - 2341.

\bibitem{CS}
{\sc J. H. Conway and N. J. A. Sloane},  {\em Sphere packings,
lattices and groups.}, 3rd ed., Springer, New York, 1999.

\bibitem{DGS}
{\sc P. Delsarte, J. M. Goethals, and J. J. Seidel}, {\em Spherical
codes and designs}, Geom. Dedicata, 6 (1977), 363--388.

\bibitem{HS}
{\sc R. H. Hardin and N. J. A. Sloane}, {\em McLaren's Improved Snub
Cube and Other New Spherical Designs in Three Dimensions}, Discrete
Comput. Geom., 15 (1996), 429--441.



\bibitem{KM}
{\sc J. Korevaar and J. L. H. Meyers}, {\em Spherical Faraday cage
for the case of equal point charges and Chebyshev-type quadrature on
the sphere}, Integral Transforms Spec. Funct. 1 (1993), 105--117.

\bibitem{SK2}
{\sc A. B. J. Kuijlaars and E. B. Saff}, {\em Asymptotics for
minimal discrete energy on the sphere}, Trans. Amer. Math. Soc., 350
(1998), 523--538.

\bibitem{MNW}
{\sc H. N. Mhaskar, F. J. Narcowich, and J. D. Ward}, {\em Spherical
Marcinkiewicz-Zygmund inequalities and positive quadrature}, Math.
Comp., 70 (2001), 1113--1130.

\bibitem{M}
{\sc O. R. Musin}, {\em The kissing number in four dimensions},
Annals of Math., 68 (2008), 1--32.


\bibitem{OCC}
{\sc Donal O'Regan, Yeol Je Cho and Yu Qing Chen}, {\em Topological
degree theory and applications}, Chapman \& Hall/CRC, 2006.

\bibitem{SK}
{\sc E. B. Saff and A. B. J. Kuijlaars}, {\em Distributing many
points on a sphere}, Math. Intelligencer, 19 (1997), 5--11.

\bibitem{SZ}
{\sc P. D. Seymour and T. Zaslavsky}, {\em Averaging sets: a
generalization of mean values and spherical designs}, Adv. Math., 52
(1984), 213--240.

\bibitem{Wag}
{\sc G. Wagner}, {\em On averaging sets}, Monatsh. Math., 111
(1991), 69--78.


\bibitem{Yu}
{\sc V. A. Yudin}, Lower bounds for spherical designs, Izv. Ross.
Akad. Nauk. Ser. Mat. 61 (1997), 211--233. English transl., Izv.
Math. 61 (1997), 673--683.

\end{thebibliography}
\end{document}